\documentclass[11pt]{article}
\usepackage{amsfonts}
\usepackage{amssymb}
\usepackage{amsmath}
\usepackage[francais]{babel}
\usepackage{graphicx}

\newtheorem{theorem}{\sc{Th\'eor\`eme}}

\newtheorem{definition}[theorem]{\sc{D\'efinition}}
\newtheorem{example}[theorem]{Exemple}

\newtheorem{proposition}[theorem]{\sc{Proposition}}

\begin{document}

\date{}
\title{\sc{Limites projectives fortes d'alg\'{e}bro\"{\i}des de Lie}
}
\author{Patrick CABAU
\footnote{L'auteur remercie Fernand Pelletier pour les diverses remarques et discussions.}\\
Laboratoire de Math\'{e}matiques \\
Universit\'{e} de Savoie, Campus scientifique \\
73376. Le Bourget du Lac\ \
France }

\maketitle

{\footnotesize

\textbf{R\'{e}sum\'{e} :} On d\'{e}finit la notion de limite projective forte d'alg\'{e}bro\"{\i}des de Lie puis on \'{e}tudie les structures de fibr\'{e}s vectoriels fr\'{e}ch\'{e}tiques associ\'{e}s et la compatibilit\'{e} entre les divers morphismes. Ce type de structure semble \^{e}tre un cadre adapt\'{e} pour diverses situations.  \\

\textbf{Abstract:} We define the notion of strong projective limit of Banach Lie
algebroids. We study the associated structures of Fr\'{e}chet bundles and the compatibility with the different morphisms. This kind of structure seems to be a convenient framework for various situations. \\

\textbf{Mots cl\'{e} :} alg\'{e}bro\"{\i}de de Lie, limite projective, fibr\'{e}
vectoriel banachique, fibr\'{e} vectoriel fr\'{e}ch\'{e}tique, calcul diff\'{e}rentiel adapt\'{e}, tenseur de Nijenhuis, distribution involutive, tenseur de Poisson, limite inverse d'espaces de Banach, semi-gerbe. \\


\textbf{Classification AMS :} 46A13, 46T05, 46T20, 58A20, 58A30, 58B10, 58B25, 26E15. \\

\quotation
}

\section{Introduction}

La notion d'alg\'{e}bro\"{\i}de de Lie a \'{e}t\'{e} introduite par J. Pradines
dans \cite{Pra} en liaison avec les groupo\"{\i}des de Lie.

Cette notion qui g\'{e}n\'{e}ralise \`{a} la fois la structure d'alg\`{e}bre
de Lie et le fibr\'{e} tangent \`{a} une vari\'{e}t\'{e} appara\^{\i}t comme
un cadre adapt\'{e} pour des probl\`{e}mes qui interviennent notamment :\\

\begin{description}
\item[--] en m\'{e}canique o\`{u} une th\'{e}orie des syst\`{e}mes lagrangiens et hamiltoniens peut \^{e}tre d\'{e}velopp\'{e}e sur de telles structures (cf. \cite{Wein2}, \cite{CLMM}, \cite{CarMar})

\item[--] en g\'{e}om\'{e}trie symplectique en vue de la symplectisation de
vari\'{e}t\'{e}s de Poisson et d'applications \`{a} la quantisation (\cite{Kar}, \cite{Wein1})

\item[--] en G\'{e}om\'{e}trie o\`{u} la classification des alg\'{e}bro\"{i}des de Lie (\cite{FerStr}) est associ\'{e}e aux G-structures de types finis ; cette notion de G-structure couvre la plupart des structures g\'{e}om\'{e}triques classiques (\cite{Mol}).

\item[--] en th\'{e}orie du contr\^{o}le optimal o\`{u} existe une version
du principe du maximum de Pontryagin (cf. \cite{Mart}).
\end{description}

\bigskip

En dimension finie, il existe une bijection entre
\begin{description}
\item[--] structures d'alg\'{e}bro\"{i}des de Lie sur un fibr\'{e} muni d'une ancre et structures de Poisson sur son dual
\item[--] structures d'alg\'{e}bro\"{i}des de Lie et diff\'{e}rentielles de Lie (cf. \cite{Marl}, \cite{CLMM}).
\end{description}
La g\'{e}n\'{e}ralisation de ces situations au contexte de vari\'{e}t\'{e}s banachiques est \'{e}tudi\'{e}e dans \cite{CabPel}.

\bigskip

L'\'{e}tude des limites projectives (ou inverses) de syst\`{e}mes  de divers types d'espaces a donn\'{e} lieu \`{a} de nombreux travaux :
\begin{description}
\item[--] limites projectives de fibr\'{e}s tangents d'une vari\'{e}t\'{e} diff\'{e}rentiable de dimension finie (cf. \cite{Gal3}) et plus g\'{e}n\'{e}ralement limites projectives de fibr\'{e}s vectoriels (cf. \cite{AghSur2}), un exemple classique \'{e}tant fourni par la g\'{e}om\'{e}trie du fibr\'{e} des jets d'ordre infini telle qu'elle est d\'{e}velopp\'{e}e par exemple dans \cite{Sau}
\item[--] limites projectives de groupes de Lie banachiques \'{e}tudi\'{e} dans \cite{Gal1} en liaison avec les groupes ILB (\cite{Omo}, \cite{Schm})
\item[--] \textit{universal laminated surfaces} \'{e}tudi\'{e}es par Nag et Sullivan dans \cite{NagSul} et utilis\'{e}es en Physique Math\'{e}matique.
\end{description}
Rappelons que la notion de limite projective a \'{e}t\'{e} introduite par Weil dans \cite{Weil} pour discuter de la structure des groupes localement compacts. \\

De nombreux probl\`{e}mes apparaissent sur des vari\'{e}t\'{e}s model\'{e}es sur des espaces de Fr\'{e}chet $\mathbb{F}$ : r\'{e}solution dans un cadre  g\'{e}n\'{e}ral d'\'{e}quations diff\'{e}rentielles (cf. \cite{Ham}) et structure pathologique du groupe $Gl(\mathbb{F})$ (qui n'admet pas de structure raisonnable de groupe de Lie). Ces probl\`{e}mes ont une solution sur certaines limites projectives d'espaces : d'une part, existence de courbes int\'{e}grales de champs de vecteurs, de courbes autoparall\`{e}les relativement \`{a} des connexions lin\'{e}aires (cf. \cite{AghSur2}), section globale horizontale pour des connexions sur certains espaces (\cite{AghSur1}) ; d'autre part, existence d'un groupe de Lie g\'{e}n\'{e}ralis\'{e} $H_{0}(\mathbb{F})$ comme groupe structural pour le fibr\'{e} tangent (cf. \cite{Gal2}).

\bigskip

On s'int\'{e}resse ici au cadre des limites projectives d'alg\'{e}bro\"{\i}des de Lie banachiques que l'on peut munir de structures fr\'{e}ch\'{e}tiques. On peut trouver dans \cite{KisVan} la notion d'alg\'ebro\"ide de Lie variationnel, utilis\'{e} dans le cadre des EDP, o\`{u} les champs de vecteurs sont remplac\'{e}s par des sections d'un fibr\'{e} au dessus d'une limite projective de jets.
\bigskip

Le r\'{e}sultat principal de ce papier (Th\'{e}or\`{e}me \ref{T_LimiteProjectiveAlgebroidesLieBanachiques}) affirme que la limite projective forte $\left( \underleftarrow{\lim }E_{i},\underleftarrow{\lim }\pi _{i},%
\underleftarrow{\lim }M_{i},\underleftarrow{\lim }\rho _{i}\right) $ d' alg\'{e}bro\"{\i}des de Lie banachiques (o\`{u} $E_i$ est un fibr\'{e} vectoriel de base $M$ et o\`{u} $\rho_i:E_i \rightarrow TM_i$ est l'ancre)  a une structure d'alg\'{e}bro\"{\i}de de Lie fr\'{e}ch\'{e}tique.
\bigskip

On rappelle dans la partie \ref{2_VarietesDimensionInfinie} les notions de
vari\'{e}t\'{e}s et de fibr\'{e}s model\'{e}s sur des espaces vectoriels
adapt\'{e}s, \textit{convenient vector spaces} selon la terminologie de Kriegel et
Michor (\cite{KriMic1}). La notion de limite projective forte de fibr\'{e}s
vectoriels banachiques introduite dans \cite{AghSur2} et g\'{e}n\'{e}ralisant
des r\'{e}sultats obtenus sur le fibr\'{e} tangent par \cite{Gal3} est rappel\'{e}e dans la partie \ref{3_LimiteProjectiveForteFVB}. La notion d'alg\'{e}bro\"{\i}de de Lie banachique est pr\'{e}sent\'{e}e ainsi que les notions de d\'{e}riv\'{e}e de Lie, de diff\'{e}rentielle ext\'{e}rieure et de morphisme
dans la partie \ref{4_AlgebroidesLieBanachiques} (cf. \cite{Ana}). On \'{e}tablit dans la partie \ref{5_LimiteProjectiveForteAlgebroidesLieBanachiques} la structure fr\'{e}ch\'{e}tique de la limite projective de ce type d'alg\'{e}bro\"{\i}des (th\'{e}or\`{e}me \ref{T_LimiteProjectiveAlgebroidesLieBanachiques}). Dans la partie \ref{6_Exemples} on donne ensuite des exemples de telles structures o\`{u}:
\begin{description}
  \item[--] $E_i=TM_i$ et l'ancre est un tenseur de Nijenhuis (cadre adapt\'{e} \`{a} l'oscillateur harmonique de dimension infinie)
  \item[--] $E_i$ est un sous-fibr\'{e} particulier de $TM_i$ :
  \begin{itemize}
    \item pour des rangs finis, on obtient la notion de diffi\'{e}t\'{e}
    \item le cadre des limites inverses d'espaces de Banach (ou de Hilbert) correspond aux rangs de dimension infinie ; on a alors un cadre int\'{e}ressant pour divers probl\`{e}mes en th\'{e}orie quantique des champs.
  \end{itemize}
\end{description}
La derni\`{e}re partie de ce travail est d\'{e}volue \`{a} l'\'{e}tude des limites projectives de semi-gerbes (eng. \textit{semisprays}) et de courbes admissibles.

\section{\label{2_VarietesDimensionInfinie}Vari\'{e}t\'{e}s de dimension
infinie model\'{e}es sur des espaces adapt\'{e}s}

Le calcul diff\'{e}rentiel classique fonctionne bien sur des vari\'{e}t\'{e}%
s de dimension finie ou banachique (cf. \cite{Lan}). Les limites projectives
de tels espaces et de fibr\'{e}s associ\'{e}s requierrent le calcul en
dimension infinie d\'{e}velopp\'{e} notamment dans \cite{KriMic1}. On
rappelle essentiellement dans cette partie\ les r\'{e}sultats \'{e}nonc\'{e}%
s dans \cite{KriMic2}, \S 2.

\subsection{Applications de classe $C^{\infty}$ sur des espaces vectoriels
adapt\'{e}s}

Si l'on souhaite munir un espace vectoriel localement convexe s\'{e}par\'{e}
$E$ d'une structure diff\'{e}rentiable, la notion premi\`{e}re est celle de
courbe $c:\mathbb{R}\rightarrow E$ de classe $C^{\infty}$.

La courbe $c$ est dite \textit{diff\'{e}rentiable} si pour tout $t\in
\mathbb{R}$, la limite du taux $\dfrac{1}{s}\left[ c\left( t+s\right)
-c\left( t\right) \right]$ existe. $c$ est dit de classe $C^{\infty}$
si toutes les d\'{e}riv\'{e}es it\'{e}r\'{e}es existent.

L'espace $\mathcal{C}=C^{\infty}\left( \mathbb{R},E\right) $ des courbes $c:%
\mathbb{R}\rightarrow E$ de classe $C^{\infty}$ ne d\'{e}pend pas de la
topologie localement convexe de $E$ mais uniquement de sa bornologie associ%
\'{e}e (syst\`{e}me de ses ensembles born\'{e}s).

L'espace $E$ est appel\'{e} espace vectoriel adapt\'{e} (\textit{convenient
vector space }selon la terminologie de \cite{KriMic1}) s'il v\'{e}rifie la
condition de $c^{\infty}$--compl\'{e}tude~:

\bigskip Une courbe $c:\mathbb{R}\rightarrow E$ est de classe $C^{\infty}$
si et seulement si $\lambda\circ c:$ $\mathbb{R}\rightarrow\mathbb{R}$ est $%
C^{\infty}$ pour tout $\lambda\in E^{\ast}$ o\`{u} $E^{\ast}$ est le dual
constitu\'{e} de toutes les applications lin\'{e}aires continues sur $E$.

La topologie finale relative \`{a} l'ensemble $\mathcal{C}$ est appel\'{e}e $%
c^{\infty}$--topologie de $E$ et est not\'{e}e $c^{\infty}E$. Un ouvert pour
cette topologie est appel\'{e} ouvert $c^{\infty}$.

Pour des espaces de Fr\'{e}chet, cette topologie co\"{\i}ncide avec la
topologie d'espace vectoriel localement convexe donn\'{e}e. Pour d'autres
espaces (e.g. l'espace $\mathcal{D}$ des fonctions tests \`{a} support
compact sur $\mathbb{R}$) elle est strictement plus fine.

\medskip

Consid\'{e}rons maintenant deux espaces vectoriels $E$ et $F$ localement
convexes et soit $U\subset E$ un ouvert $c^{\infty}$. Une application $%
f:E\supset U\rightarrow F$, o\`{u} $E$ et $F$ sont deux espaces vectoriels
adapt\'{e}s et o\`{u} $U$ est un $c^{\infty}$--ouvert de $E$, est dite de
classe $C^{\infty}$ si $f\circ c\in C^{\infty}\left( \mathbb{R},F\right) $
pour toute $c\in C^{\infty}\left( \mathbb{R},U\right) .$ De plus, cf. \cite%
{KriMic2}, 2.3 (5), l'espace $C^{\infty}\left( U,F\right) $ peut \^{e}tre
aussi muni d'une structure d'espace vectoriel adapt\'{e}.

\subsection{Vari\'{e}t\'{e}s diff\'{e}rentiables}

\subsubsection{\label{D_VarieteDifferentiableConvenient}Structure de vari%
\'{e}t\'{e} diff\'{e}rentiable sur un ensemble}

Une carte $\left( U,\varphi\right) $ sur un ensemble $M$ est une bijection $%
\varphi:U\rightarrow\varphi\left( U\right) \subset E$ d'une partie $U$ de $M$
sur un $c^{\infty}$--ouvert d'un espace vectoriel adapt\'{e} $E$.

Une famille $\left( U_{\alpha},\varphi_{\alpha}\right) _{\alpha\in A}$ de
cartes est appel\'{e}e un \textit{atlas} si tous les changements de cartes $%
\varphi_{\alpha\beta}=\varphi_{\alpha}\circ\left( \varphi_{\beta}\right)
^{-1}:\varphi_{\beta}\left( U_{\alpha}\cap U_{\beta}\right) \rightarrow
\varphi_{\alpha}\left( U_{\alpha}\cap U_{\beta}\right) $ sont de classe $%
C^{\infty}$.

Deux atlas sont dits \'{e}quivalents si leur r\'{e}union est encore un atlas.

L'ensemble $M$ muni d'une classe d'\'{e}quivalence d'atlas est appel\'{e}
\textit{vari\'{e}t\'{e} diff\'{e}rentiable} de classe $C^{\infty}$.

Un sous-ensemble $W$ de la vari\'{e}t\'{e} $M$ est ouvert si et seulement si
pour tout $\alpha\in A$ le sous ensemble $\varphi_{\alpha}\left( U_{\alpha
}\cap W\right) $ de $E$ est $c^{\infty}$--ouvert.

La topologie ainsi d\'{e}finie est alors la topologie finale relativement
\`{a} l'espace des courbes de classe $C^{\infty }$.

\subsubsection{Sous-vari\'{e}t\'{e}s}

Un sous-ensemble $N$ d'une vari\'{e}t\'{e} $M$ est appel\'{e}e \textit{%
sous-vari\'{e}t\'{e}} si pour tout $x\in N$ il existe une carte $\left(
U,\varphi \right) $ de $M$ telle que%
\begin{equation*}
\varphi \left( U\cap N\right) =\varphi \left( U\right) \cap F
\end{equation*}

o\`{u} $F$ est un sous-espace vectoriel ferm\'{e} de l'espace vectoriel adapt%
\'{e} $E$.

\subsubsection{\label{D_AppClasseCInfini}Applications de classe $C^{\infty}$
entre vari\'{e}t\'{e}s}

Une application $f:M\rightarrow N$ entre deux vari\'{e}t\'{e}s diff\'{e}%
rentiables est dite de classe $C^{\infty}$ si pour tout $x\in M$ et pour
toute carte $\left( V,\psi\right) $ de $N$ telle que $f\left( x\right) \in V$
il existe une carte $\left( U,\varphi\right) $ de $M$ telle que $x\in U,$ $%
f\left( U\right) \subset V$ et telle que $\psi\circ f\circ\varphi^{-1}$ est
de classe $C^{\infty}$. Ceci est \'{e}quivalent au fait que $f\circ c$ est $%
C^{\infty}$ pour chaque courbe $c:\mathbb{R}\rightarrow M$ de classe $%
C^{\infty}$.

On note $\mathcal{F}$ l'anneau des fonctions de classe $C^{\infty}$ de $M$
dans $\mathbb{R}$.

\subsubsection{Fibr\'{e}s vectoriels}

Soit $p:F\rightarrow M$ une application $C^{\infty}$ entre vari\'{e}t\'{e}s
diff\'{e}rentiables $F$ et $M$.

Une \textit{carte de fibr\'{e} vectoriel} sur $\left( F,p,M\right) $ est un
couple $\left( U,\Phi\right) $ o\`{u} $U$ est un ouvert de $M$ et o\`{u} $%
\Phi$ est un diff\'{e}omorphisme respectant la fibre, i.e. pour lequel le
diagramme ci-dessous est commutatif%
\begin{equation*}
\begin{array}{ccccc}
F_{|U}=p^{-1}\left( U\right) &  & \overset{\Phi}{\longrightarrow} &  &
U\times V \\
& p\searrow &  & \swarrow pr1 &  \\
&  & U &  &
\end{array}%
\end{equation*}
o\`{u} $V$ est un espace vectoriel adapt\'{e} fixe appel\'{e} fibre
strandard.

Deux cartes $\left( U_{1},\Phi_{1}\right) $ et $\left( U_{2},\Phi
_{2}\right) $ sont dites \textit{compatibles} si $\Phi_{1}\circ\left(
\Phi_{2}\right) ^{-1}$ $\left( x,v\right) $ peut \^{e}tre \'{e}crit sous la
forme $\left( x,\Phi_{1,2}\left( x\right) \left( v\right) \right) $ o\`{u} $%
\Phi_{1,2}:U_{1}\cap U_{2}\rightarrow GL\left( V\right) $ . La fonction $%
\Phi_{1,2}$ est alors unique et $C^{\infty}$ dans $L\left( V\right) $ o\`{u}
$L\left( V\right) $ est l'espace des applications lin\'{e}aires born\'{e}es
(donc $C^{\infty}$)$.$

Un \textit{atlas de fibr\'{e} vectoriel} $\left( U_{\alpha},\Phi_{\alpha
}\right) _{\alpha\in A}$ pour $p:F\rightarrow M$ est un ensemble de cartes $%
\left( U_{\alpha},\Phi_{\alpha}\right) $ deux \`{a} deux compatibles o\`{u} $%
\left( U_{\alpha}\right) _{\alpha\in A}$ est un recouvrement ouvert de la
vari\'{e}t\'{e} $M$. Deux atlas sont \textit{\'{e}quivalents} si leur r\'{e}%
union est encore un atlas.

Un \textit{fibr\'{e} vectoriel}, de classe $C^{\infty}$, $p:F\rightarrow M$
est la donn\'{e}e de vari\'{e}t\'{e}s $F$ (espace total), $M$ (base) et
d'une application $p:F\rightarrow M$ (projection) de classe $C^{\infty}$
munies d'une classe d'\'{e}quivalence d'atlas.

Une section $s$ de $p:F\rightarrow M$ est une application $C^{\infty}$ $%
u:M\rightarrow F$ tellle que $p\circ u=\lim{Id}_{M}$.

L'espace \underline{$F$} des sections de $F$ peut \^{e}tre muni d'une
structure d'espace vectoriel adapt\'{e}.

\subsubsection{Vecteurs tangents}

Un \textit{vecteur tangent (cin\'{e}matique)} en un point $x$ d'une vari\'{e}%
t\'{e} $M$ est une classe d'\'{e}quivalence pour la relation suivante%
\begin{equation*}
c_{1}\sim c_{2}\text{ si et seulement si }\left\{
\begin{array}{c}
c_{1}\left( 0\right) =c_{2}\left( 0\right) =x\in U \\
\left( \varphi \circ c_{1}\right) ^{\prime }\left( 0\right) =\left( \varphi
\circ c_{2}\right) ^{\prime }\left( 0\right)%
\end{array}%
\right.
\end{equation*}

o\`{u} $\left( U,\varphi \right) $ est une carte de $M$ centr\'{e}e en $x$.

Il existe une autre notion de vecteur tangent qui est la notion de \textit{%
vecteur tangent op\'{e}rationnel }(\cite{KriMic1}, 28.1) qui ne co\"{\i}%
ncide pas n\'{e}cessairement sur des vari\'{e}t\'{e}s adapt\'{e}es \`{a} la
notion de vecteur tangent cin\'{e}matique.

\subsubsection{Fibr\'{e} tangent}

L'ensemble de tous les vecteurs tangents en les divers points de la vari\'{e}%
t\'{e} $M$ peut \^{e}tre muni d'une structure de fibr\'{e} vectoriel~; il
est alors appel\'{e} \textit{fibr\'{e} tangent (cin\'{e}matique)} \`{a} $M$
et not\'{e} $TM$.

Si $\left( U_{\alpha },\varphi _{\alpha }\right) _{\alpha \in A}$ est un
atlas de la vari\'{e}t\'{e} $M$ alors les changements de cartes dans $TM$
font intervenir les diff\'{e}rentielles $d\varphi _{\alpha \beta }$.

\subsubsection{Champs de vecteurs}

Un \textit{champ de vecteurs cin\'{e}matique} est une section $C^{\infty}$
du fibr\'{e} tangent cin\'{e}matique $TM$. On note $\mathfrak{X}\left(
M\right) $ l'espace des champs de vecteurs cin\'{e}matiques qui peut lui
aussi \^{e}tre muni d'une structure d'espace vectoriel adapt\'{e}.

Pour des vari\'{e}t\'{e}s r\'{e}guli\`{e}res au sens de \cite{KriMic1}, 14,
le crochet de Lie de deux \'{e}l\'{e}ments $X$ et $Y$ de $\mathfrak{X}\left(
M\right) $ peut \^{e}tre d\'{e}fini, en supposant que $M$ est un ouvert $%
c^{\infty}$ d'un espace vectoriel adapt\'{e} $E$, par%
\begin{equation*}
\left[ X,Y\right] =dY\left( X\right) -dX\left( Y\right)
\end{equation*}

o\`{u} $X$ et $Y$ sont alors vus comme application $M\rightarrow E$ de classe
$C^{\infty}$.

\subsubsection{Application tangente}

Si $M$ et $N$ sont deux vari\'{e}t\'{e}s diff\'{e}rentiables, \`{a} toute
application $f:M\rightarrow N$ de classe $C^{\infty }$, on peut associer
\`{a} tout point $x\in M$ une application lin\'{e}aire $T_{x}f:T_{x}M%
\rightarrow T_{f\left( x\right) }M$ qui associe au vecteur tangent \`{a} une
courbe $c$ passant par $x$ le vecteur tangent \`{a} la courbe $f\circ c$
passant par $f\left( x\right) $.

L'application obtenue $Tf:TM\rightarrow TN$ est alors $C^{\infty }$ et est
appel\'{e}e \textit{application tangente} de $f$.

\subsubsection{Champs de vecteurs reli\'{e}s}

$\bigskip $Soient $M$ et $N$ deux vari\'{e}t\'{e}s diff\'{e}rentiables et $%
f:M\rightarrow N$ une application $C^{\infty }$. On dit que deux champs de
vecteurs cin\'{e}matiques sont $f$--reli\'{e}s si%
\begin{equation*}
Tf\circ X=Y\circ f
\end{equation*}

$X_{1},X_{2},Y_{1},Y_{2}$ d\'{e}signent maintenant 4 champs de vecteurs.

Soit $f:M\rightarrow N$ une application de classe $C^{\infty}$. Si $X_{1}$
et $Y_{1}$ (resp. $X_{2}$ et $Y_{2}$) sont $f$--reli\'{e}s alors $\left[
X_{1},X_{2}\right] $ et $\left[ Y_{1},Y_{2}\right] $ sont aussi $f$--reli%
\'{e}s.

Si $f:M\rightarrow N$ est un diff\'{e}omorphisme local alors pour tout champ
de vecteurs $Y\in \mathfrak{X}\left( N\right) $ il existe un champ de
vecteurs $f^{\ast }Y\in \mathfrak{X}\left( M\right) $ d\'{e}fini par $\left(
f^{\ast }Y\right) \left( x\right) =\left( T_{x}f\right) ^{-1}\left( Y\left(
f\left( x\right) \right) \right) .$ L'application lin\'{e}aire $f^{\ast }:%
\mathfrak{X}\left( N\right) \rightarrow \mathfrak{X}\left( M\right) $ est
alors un homomorphisme d'alg\`{e}bres de Lie.

\subsubsection{Courbe int\'{e}grale d'un champ de vecteurs cin\'{e}matique}

Une courbe $c:I\rightarrow M$ de classe $C^{\infty}$ est dite \textit{courbe
int\'{e}grale} du champ de vecteurs $X$ si pour tout $t\in I,$ on a~:~$%
c^{\prime}\left( t\right) =X\left( c\left( t\right) \right) $.

Pour un champ de vecteurs donn\'{e} les courbes int\'{e}grales peuvent ne
pas exister (\cite{KriMic1}, 32.12, exemple 1) localement et m\^{e}me si
elles existent, elles peuvent de pas \^{e}tre uniques pour une condition
initiale donn\'{e}e (\cite{KriMic1}, 32.12, exemple 2). Ceci est d\^{u} au
fait que les r\'{e}sultats classiques relatifs \`{a} l'existence et \`{a}
l'unicit\'{e} de solutions d'\'{e}quations diff\'{e}rentielles sont fond\'{e}%
s sur des th\'{e}or\`{e}mes r\'{e}sultant du th\'{e}or\`{e}me du point fixe
sur des espaces de Banach.

\subsubsection{Flot d'un champ de vecteurs}

Un \textit{flot local pour un champ de vecteurs} $X\in\mathfrak{X}\left(
M\right) $ est une application, de classe $C^{\infty}$, $\varphi^{X}:M\times%
\mathbb{R}\supset U\rightarrow M$ d\'{e}finie sur un ouvert $c^{\infty}$ $U$
de $M\times\left\{ 0\right\} $ tel que~:

\begin{enumerate}
\item $U\cap\left( \left\{ x\right\} \times\mathbb{R}\right) $ est un
intervalle ouvert connexe

\item Si $\varphi^{X}\left( x,s\right) $ existe alors $\varphi^{X}\left(
x,t+s\right) $ existe si et seulement si $\varphi^{X}\left(
\varphi^{X}\left( x,s\right) ,t\right) $
existe et on a $\varphi^{X}\left( x,t+s\right) =\varphi^{X}\left(
\varphi^{X}\left( x,s\right) ,t\right) $

\item $\varphi^{X}\left( x,0\right) =x$ pour tout $x\in M$

\item $\dfrac{d}{dt}\varphi^{X}\left( x,t\right) =X\left( \varphi ^{X}\left(
x,t\right) \right) $
\end{enumerate}

Si un champ de vecteurs cin\'{e}matique $X$ admet un flot local $\varphi
^{X} $ alors pour toute courbe int\'{e}grale $c$ de $X$, on a $c\left(
t\right) =\varphi ^{X}\left( c\left( 0\right) ,t\right) $ et il existe ainsi
un unique flot maximal.

\subsubsection{Fibr\'{e} cotangent}

Une $1-$forme en un point $x$ d'une vari\'{e}t\'{e} $M$ est une forme lin%
\'{e}aire born\'{e}e sur l'espace vectoriel adapt\'{e} $T_{x}M$ (donc
appartenant \`{a} $T_{x}M^{\prime }$). L'ensemble de toutes ces $1-$formes en
les divers points de $M$ peut \^{e}tre muni d'une structure de fibr\'{e}
vectoriel appel\'{e} \textit{fibr\'{e} cotangent (cin\'{e}matique)} et not\'{e} $T^{\prime }M$.

Un atlas $\left( U_{\alpha },\varphi _{\alpha }\right) _{\alpha \in A}$ de
classe $C^{\infty }$ de $M$ donne naissance aux fonctions de transition $%
x\mapsto d\left( \varphi _{\beta }\circ \left( \varphi _{\alpha }\right)
^{-1}\right) _{\varphi _{\alpha }\left( x\right) }$.

\subsubsection{Formes diff\'{e}rentielles}

\bigskip Sur une vari\'{e}t\'{e} $M$ une $1-$\textit{forme diff\'{e}%
rentielle (cin\'{e}matique)} n'est autre qu'une section $C^{\infty }$ du fibr\'{e} cotangent cin\'{e}matique $T^{\prime }M$. L'ensemble de ces formes diff\'{e}rentielles peut \^{e}tre muni d'une structure d'espace vectoriel adapt\'{e}.

\bigskip

Sur une vari\'{e}t\'{e} $C^{\infty }$ r\'{e}guli\`{e}re $M$, la classe de
formes diff\'{e}rentielles (\cite{KriMic1}, 33.22) stable par d\'{e}riv\'{e}%
e de Lie $L_{X}$, d\'{e}rivation ext\'{e}rieure $d$, produit int\'{e}rieur $%
i_{X}$ et image r\'{e}ciproque $f^{\ast }$ est%
\begin{equation*}
\Omega ^{k}\left( M\right) =\underline{L_{\text{alt}}^{k}\left( TM,M\times
\mathbb{R}\right) }
\end{equation*}

La d\'{e}riv\'{e}e de Lie $L:\mathfrak{X}\left( M\right) \times\Omega
^{k}\left( M\right) \rightarrow\Omega^{k}\left( M\right) $ est une
application $C^{\infty}$ d\'{e}finie par%
\begin{equation*}
\left( L_{X}\omega\right) \left( X_{1},\dots,X_{k}\right) =X\left(
\omega\left( X_{1},\dots,X_{k}\right) \right) -{\displaystyle%
\sum\limits_{i=1}^{k}} \omega\left( X_{1},\dots,\left[ X,X_{i}\right]
,\dots,X_{k}\right)
\end{equation*}

La diff\'{e}rentielle ext\'{e}rieure $d:\Omega^{k}\left( M\right)
\rightarrow\Omega^{k+1}\left( M\right) $ est $C^{\infty}$ et est d\'{e}finie
par%
\begin{align*}
\left( d\omega\right) \left( x\right) \left( X_{0},\dots,X_{k}\right) & = {%
\displaystyle\sum\limits_{i=0}^{k}} \left( -1\right) ^{i}X_{i}\left(
\omega\left( X_{0},\dots,\widehat{X_{i}},\dots,X_{k}\right) \right) \\
& + {\displaystyle\sum\limits_{0\leq i<j\leq k}} \left( -1\right)
^{i+j}\omega\left( \left[ X_{i},X_{j}\right] ,X_{0},\dots,\widehat{X_{i}}%
,\dots,\widehat{X_{j}},\dots,X_{k}\right)
\end{align*}

\section{\label{3_LimiteProjectiveForteFVB}Limite projective forte de fibr\'{e}s vectoriels banachiques}

\subsection{Limites projectives d'espaces topologiques}

Un syst\`{e}me projectif d'espaces topologiques est une suite $\left( \left(
X_{i},\delta_{i}^{j}\right) _{j\geq i}\right) _{i\in\mathbb{N}}$ o\`{u}

\begin{description}
\item[--] pour tout $i\in\mathbb{N},$ $X_{i}$ est un espace topologique

\item[--] pour tous $i,$ $j$ $\in\mathbb{N},$ tels que $j\geq i,$ $\delta
_{i}^{j}:X_{j}\rightarrow X_{i}$ est une application continue

\item[--] pour tout $i\in \mathbb{N},$ $\delta _{i}^{i}={Id}_{X_{i}}$

\item[--] pour tous entiers naturels $i\leq j\leq k$, $\delta_{i}^{j}\circ%
\delta_{j}^{k}=\delta_{i}^{k}$.
\end{description}

Un \'{e}l\'{e}ment $\left( x_{i}\right) _{i\in\mathbb{N}}$ du produit ${%
\displaystyle\prod\limits_{i\in\mathbb{N}}} X_{i}$ est appel\'{e} un fil
[thread] si pour tous $j\geq i$, $\delta_{i}^{j}\left( x_{j}\right) =x_{i}$.

Le sous-espace X de $\displaystyle\prod\limits_{i\in\mathbb{N}}X_{i}$ constitu\'{e} de tels \'{e}l\'{e}ments, muni de la topologie projective, i.e. de la topologie la moins fine rendant continues toutes les applications $\delta_{i}=p_{i/X}$ (o\`{u} $p_{i}: \displaystyle\prod\limits_{k\in\mathbb{N}}X_{k}\rightarrow X_{i}$ d\'{e}signe la projection sur $X_{i}$) est appel\'{e}e limite projective de la suite $((X_{i},\delta_{i}^{j})_{j\geq i})_{i\in\mathbb{N}}$ (cf. \cite{Scha}, p.52). On note $X=\underleftarrow{\lim}X_{i}$ et $\cal{T}$ la topologie projective.
Si $x\in X$ et si $x_{i}=\delta_{i}(x)$, une base de voisinages de $\cal{T}$ est alors donn\'{e}e par toutes les intersections
 ${\displaystyle\bigcap\limits_{l\in L}}\left(  \delta_{l}\right)  ^{-1}\left(  U_{l}\right)$
o\`{u} $U_{l}$ est un voisinage de $x_{l}$ relativement \`{a} la topologie de l'espace $X_l$ et o\`{u} $L$ est un ensemble fini.

Soient $\left( \left( X_{i},\delta_{i}^{j}\right) _{j\geq i}\right) _{i\in%
\mathbb{N}}$ et $\left( \left( Y_{i},\gamma_{i}^{j}\right) _{j\geq i}\right)
_{i\in\mathbb{N}}$ deux syst\`{e}mes projectifs de limites projectives
respectives $X$ et $Y$.

Une suite $\left( f_{i}\right) _{i\in\mathbb{N}}$ d'applications continues $%
f_{i}:X_{i}\rightarrow Y_{i}$ v\'{e}rifiant, pour tous $i,j\in\mathbb{N},$ $%
j\geq i,$ la condition de coh\'{e}rence%
\begin{equation*}
\gamma_{i}^{j}\circ f_{j}=f_{i}\circ\delta_{i}^{j}
\end{equation*}

est appel\'{e}e \textit{suite projective d'applications}.

La limite projective de cette suite est l'application%
\begin{equation*}
\begin{array}{cccc}
f: & X & \rightarrow & Y \\
& \left( x_{i}\right) _{i\in\mathbb{N}} & \mapsto & \left( f_{i}\left(
x_{i}\right) \right) _{i\in\mathbb{N}}%
\end{array}%
\end{equation*}

L'application $f$ est alors continue et est un hom\'{e}omorphisme si tous
les $f_{i}$ sont eux-m\^{e}mes des hom\'{e}omorphismes (cf. \cite{AbbMan}).

\subsection{Limite projective forte de vari\'{e}t\'{e}s banachiques}

La suite $\left( \left( M_{i},\delta_{i}^{j}\right) _{j\geq i}\right) _{i\in
\mathbb{N}}$ est appel\'{e}e \textit{syst\`{e}me projectif fort de vari\'{e}t\'{e}s banachiques} si

\begin{description}
\item[--] $M_{i}$ est une vari\'{e}t\'{e} diff\'{e}rentiable model\'{e}e sur
l'espace de Banach $\mathbb{M}_{i}$

\item[--] $\left( \left( \mathbb{M}_{i},\overline{\delta_{i}^{j}}\right)
_{j\geq i}\right) _{i\in\mathbb{N}}$ est un syst\`{e}me projectif d'espaces
banachiques

\item[--] pour tout $x=\left( x_{i}\right) \in M=\underleftarrow{\lim }%
M_{i}, $ il existe un syst\`{e}me projectif de cartes locales $\left(
U_{i},\varphi _{i}\right) _{i\in \mathbb{N}}$ tel que $x_{i}\in U_{i}$ o\`{u}
la relation de coh\'{e}rence $\varphi _{i}\circ \delta _{i}^{j}=\overline{%
\delta _{i}^{j}}\circ \varphi _{j}$ est v\'{e}rifi\'{e}e

\item[--] $U=\underleftarrow{\lim }U_{i}$ est ouvert dans $M$.
\end{description}

La limite projective $M=\underleftarrow{\lim }M_{i}$ a alors une structure
de vari\'{e}t\'{e} fr\'{e}ch\'{e}tique model\'{e}e sur l'espace de Fr\'{e}%
chet $\mathbb{M}=\underleftarrow{\lim }\mathbb{M}_{i}$ o\`{u} la structure
diff\'{e}rentiable est d\'{e}finie via les cartes $\left( U,\varphi \right) $
o\`{u} $\varphi =\underleftarrow{\lim }\varphi _{i}:U\rightarrow \left(
\varphi _{i}\left( U_{i}\right) \right) .$\newline
$\varphi $ est bien un hom\'{e}omorphisme (limite projective d'hom\'{e}%
omorphismes) et les applications de changements de cartes
$$\left( \psi \circ
\varphi ^{-1}\right) _{|\varphi \left( U\right) }=\underleftarrow{\lim }\left( \left( \psi _{i}\circ \left( \varphi
_{i}\right) ^{-1}\right) _{|\varphi _{i}\left( U_{i}\right) }\right) $$ entre
ouverts d'espaces de Fr\'{e}chet sont de classe $C^{\infty }$ au sens de
Kriegel et Michor (cf. \cite{KriMic1}).

\begin{example}
Espace des jets d'ordre infini des sections d'un fibr\'{e} vectoriel de rang
fini au-dessus d'une vari\'{e}t\'{e} r\'{e}elle de dimension finie (cf \cite%
{Sau}, \cite{AbbMan}).\newline
\end{example}

\begin{example}
Limite projective de groupes de Banach-Lie (cf. \cite{Gal1}, \cite{Omo},
\cite{AbbMan}).\newline
Un groupe $G$ est appel\'{e} limite projective de groupes de Banach-Lie model%
\'{e} sur la limite projective $\mathbb{G}=\underleftarrow{\lim}\mathbb{G}%
_{i}$ si

\begin{enumerate}
\item $G=\underleftarrow{\lim }G_{i}$ o\`{u} $\left( G_{i},\delta
_{i}^{j}\right)$ est un syst\`{e}me projectif de groupes de Banach-Lie o\`{u}
$G_{i}$ est model\'{e} sur $\mathbb{G}_{i}$

\item Pour tout $i \in\mathbb{N}$ il existe une carte $\left(
U_{i},\varphi_{i}\right) $ centr\'{e}e en le neutre $e_{i} \in G_{i}$ telle
que~:

\begin{enumerate}
\item $\delta_{i}^{j}(U_{j})\subset U_{i}$ pour $j\geq i$

\item $\overline{\delta_{i}^{j}} \circ\varphi_{j} = \varphi_{j} \circ
\delta_{i}^{j}$

\item $\underleftarrow{\lim}\varphi_{i}(U_{i})$ est ouvert dans $\mathbb{G}$
et $\underleftarrow{\lim} U_{i}$ est ouvert dans $G$ relativement \`{a} la
topologie de limite projective.
\end{enumerate}
\end{enumerate}

\begin{itemize}
\item[--] A titre d'exemple simple, l'espace des suites r\'{e}elles $\mathbb{%
R}^{\mathbb{N}}$ muni de la topologie produit est un groupe de Lie ab\'{e}%
lien, limite projective des groupes de Lie ab\'{e}liens $\mathbb{R}^{j}$, $j
\in\mathbb{N}$.

\item[--] Des exemples plus cons\'{e}quents correspondent aux groupes de Lie
compacts eu \'{e}gard au fait que tout groupe de Lie compact est la limite
projective d'une famille de groupes de Lie compacts (cf. \cite{Weil} )
\end{itemize}

On peut notamment d\'{e}finir sur de tels groupes de Fr\'{e}chet-Lie $G$
l'exponentielle $exp_{G}$ comme limite projective de la suite $exp_{G_{i}}$.
Cette application est alors continue.
\end{example}

\subsection{Limite projective forte de fibr\'{e}s vectoriels}

Soit $\left( \left( M_{i},\delta_{i}^{j}\right) _{j\geq i}\right) _{i\in%
\mathbb{N}}$ un syst\`{e}me projectif fort de vari\'{e}t\'{e}s banachiques o%
\`{u} chaque vari\'{e}t\'{e} $M_{i}$ est model\'{e}e sur l'espace de Banach $%
\mathbb{M}_{i}$.

On consid\`{e}re pour tout entier $i$ le fibr\'{e} vectoriel banachique $%
\left( E_{i},\pi_{i},M_{i}\right) $ de fibre type $\mathbb{E}_{i}$ o\`{u},
de plus, $\left( \left( \mathbb{E}_{i},\lambda_{i}^{j}\right) _{j\geq
i}\right) _{i\in\mathbb{N}}$ constitue un syst\`{e}me projectif d'espace de
Banach.

Le syst\`{e}me $\left( \left( E_{i},f_{i}^{j}\right) _{j\geq i}\right)
_{i\in \mathbb{N}}$ est appel\'{e} \textit{syst\`{e}me projectif fort de fibr%
\'{e}s vectoriels banachiques} sur $\left( \left( M_{i},\delta
_{i}^{j}\right) _{j\geq i}\right) $ si pour tout $\left( x_{i}\right) $ il
existe un syst\`{e}me projectif de trivialisations $\left( U_{i},\tau
_{i}\right) $ de $\left( E_{i},\pi _{i},M_{i}\right) ,$ o\`{u} $\tau
_{i}:\left( \pi _{i}\right) ^{-1}\left( U_{i}\right) \rightarrow U_{i}\times
\mathbb{E}_{i}$ sont des diff\'{e}omorphismes locaux, tel que $x_{i}\in
U_{i} $ (ouvert de $M_{i}$) et o\`{u} $U=\underleftarrow{\lim }U_{i}$ est
ouvert dans $M$ avec pour tous $i,j\in \mathbb{N}$ tels que $j\geq i$ on ait
la condition de coh\'{e}rence%
\begin{equation*}
\left( \delta _{i}^{j}\times \lambda _{i}^{j}\right) \circ \tau _{j}=\tau
_{i}\circ f_{i}^{j}
\end{equation*}

On a alors la proposition suivante~qui g\'{e}n\'{e}ralise le r\'{e}sultat de
\cite{Gal3} relatif \`{a} la limite projective de fibr\'{e}s tangents \`{a}
des vari\'{e}t\'{e}s banachiques et dont on trouvera la d\'{e}monstration
dans \cite{AghSur2}.

\begin{proposition}
\label{P_LimProjFVB}Soit $\left( E_{i},\pi_{i},M_{i}\right) _{i\in \mathbb{N}%
}$ un syst\`{e}me projectif fort de fibr\'{e}s vectoriels banachiques.
\newline
Alors $\left( \underleftarrow{\lim}E_{i},\underleftarrow{\lim}\pi_{i},%
\underleftarrow{\lim}M_{i}\right) $ est un fibr\'{e} vectoriel fr\'{e}ch\'{e}%
tique.
\end{proposition}

Notons que le groupe lin\'{e}aire $Gl\left( \mathbb{E}\right) $ ne pouvant
\^{e}tre muni d'une structure de groupe de Lie ne peut pas jouer le r\^{o}le
de groupe structural. Il est ici remplac\'{e} par le groupe de Lie g\'{e}n%
\'{e}ralis\'{e}, $H^{0}\left( \mathbb{E} \right) =\underleftarrow {\lim}%
H_{i}^{0}\left( \mathbb{E}\right) $ limite projective des groupes de Lie
banachiques
\begin{equation*}
H_{i}^{0}\left( \mathbb{E}\right) =\left\{ \left( h_{1},\dots ,h_{i}\right)
\in{\displaystyle\prod\limits_{j=1}^{i}} Gl\left( \mathbb{E}_{j}\right)
:\lambda_{k}^{j}\circ h_{j}=h_{k}\circ\lambda_{k}^{j},\text{ pour }k\leq
j\leq i\right\}
\end{equation*}

On obtient alors la diff\'{e}rentiabilit\'{e} des fonctions de transition $%
\mathtt{T}$.

\section{\label{4_AlgebroidesLieBanachiques}Alg\'{e}bro\"{\i}des de Lie
banachiques}

\subsection{\label{DE_AlgebroidesLieBanachiques}D\'{e}finition. Exemples}

Soit $\pi:E\rightarrow M$ un fibr\'{e} vectoriel banachique dont la fibre
type est un espace de Banach $\mathbb{E}$.

Un morphisme de fibr\'{e}s vectoriels $\rho:E\rightarrow TM$ est appel\'{e}
\textit{ancre}. Ce morphisme induit une application \underline{$\rho$}$:%
\underline{E}\rightarrow$\underline{$TM$}=$\mathfrak{X}(M)$ d\'{e}finie pour
tout $x$ de $M$ et toute section $s$ de $E$ par : $\left( \underline {\rho}%
\left( s\right) \right) \left( x\right) =\rho\left( s\left( x\right) \right)
.$

On suppose qu'il existe un crochet $\left[ .,.\right] _{E}$ sur l'espace
\underline{$E$} qui le munisse d'une structure d'alg\`{e}bre de Lie.

\begin{definition}
Le quadruplet $\left( E,\pi,M,\rho\right) $ est appel\'{e} alg\'{e}bro\"{\i}%
de de Lie banachique si :

\begin{enumerate}
\item $\underline{\rho}$ :$\left( \underline{E},\left[ .,.\right]
_{E}\right) \rightarrow\left( \mathfrak{X}(M) ,\left[ .,.\right] \right) $
est un homomorphisme d'alg\`{e}bres de Lie

\item $\left[ s_{1},fs_{2}\right] _{E}=f\left[ s_{1},s_{2}\right]
_{E}+\left( \underline{\rho}\left( s_{1}\right) \right) \left( f\right) \
s_{2}$ pour tous $f\in\mathcal{F}$ et s$_{1},s_{2}\in$\underline{$E$}
\end{enumerate}
\end{definition}

\begin{example}
{\label{EX_AlgebroideLie_TenseurNijenhuis} $E=TM$ et $\rho=N$ est un tenseur
de Nijenhuis, i.e. v\'{e}rifiant la propri\'{e}t\'{e}\newline
\begin{equation*}
\left[ NX,NY \right] = N\left( \left[ NX,Y \right] + \left[ X,NY \right] - N
\left( \left[ N,Y \right] \right) \right)
\end{equation*}
$\left( TM,\pi,M,N \right) $ est un alg\'{e}bro\"{\i}de de Lie pour le
crochet $\left[ .,. \right] _{N}$ d\'{e}fini par \newline
\begin{equation*}
\left[ X,Y \right] _{N} = \left[ NX,Y \right] + \left[ X,NY \right] - N
\left( \left[ X,Y \right] \right)
\end{equation*}
Le cas trivial correspond au cas o\`{u} $N=\lim{Id}_{TM}$ }
\end{example}

\begin{example}
{\label{EX_AlgebroideLie_DistributionInvolutive} $E$ est une distribution,
i.e. un sous-fibr\'{e} vectoriel de $TM$ dont la fibre $E_{x}$ au dessus
d'un point $x$ de la base $M$ est un sous-espace vectoriel banachique ferm%
\'{e} de codimension finie, qui est de plus involutive. L'ancre est alors
l'injection canonique $\rho:E\rightarrow TM$. }
\end{example}

\begin{example}
{\label{EX_AlgebroideLie_TenseurPoisson} $E$ est le fibr\'{e} cotangent \`{a}
une vari\'{e}t\'{e} banachique et $\rho=P$ est un tenseur de Poisson. Le
crochet sur les sections de $T^{\ast}M$ est d\'{e}fini (cf. \cite{MagMor}) par
\[
\left\{  \alpha,\beta\right\}  _{P}=L_{P\beta}\left(  \alpha\right)
-L_{P\alpha}\left(  \beta\right)  +d\left\langle \beta,P\alpha\right\rangle
\]
 Le fait que $\left( T^{\ast }M,P,M,\left\{
.,.\right\} _{P}\right) $ ait une structure d'alg\'{e}bro\"{\i}de de Lie r%
\'{e}sulte de la propri\'{e}t\'{e}%
\begin{equation*}
\left\{ \alpha,f.\beta\right\} _{P}=f.\left\{ \alpha,\beta\right\}
_{P}+L_{P\alpha}\left( f\right) .\beta
\end{equation*}
On peut trouver une g\'{e}n\'{e}ralisation aux structures de Jacobi,
structures introduites par Lichnerowicz (\cite{Lic}), dans \cite{Pon}. }
\end{example}

\begin{example}
{\label{EX_AlgebroideLie_ActionGroupeLie} Soit une action \`{a} droite $%
\psi:M\times G\rightarrow M$ d'un groupe de Lie $G$ sur une vari\'{e}t\'{e}
banachique $M.$ On note $\mathcal{G}$ l'alg\`{e}bre de Lie de $G$. Il existe
alors un morphisme naturel du fibr\'{e} banachique trivial $M\times \mathcal{%
G}$ dans $M$ d\'{e}fini par%
\begin{equation*}
\Psi(x,X)=T_{(x,e)}\psi(0,X)
\end{equation*}
\newline
\newline
Pour tous $X$ et $Y$ dans $\mathcal{G}$, on a~:%
\begin{equation*}
\Psi(\{X,Y\})=[\Psi(X),\Psi(Y)]
\end{equation*}
\newline
\newline
o\`{u} $\{\;,\;\}$ d\'{e}signe le crochet de Lie sur $\mathcal{G}$ (cf. \cite%
{KriMic1}, 36.12).\newline
\newline
$(M\times \mathcal{G},\Psi,M,\{\;,\;\})$ a alors une structure d'alg\`{e}bro%
\"{\i}de de Lie. }
\end{example}

\subsection{Op\'{e}rateurs de d\'{e}rivation}

Peuvent \^{e}tre d\'{e}finies, sur un alg\'{e}bro\"{\i}de de Lie banachique
les notions de d\'{e}riv\'{e}e de Lie $L_{s}^{\rho}$ relativement \`{a} une
section $s$ de $E$ (cette notion g\'{e}n\'{e}ralisant la notion de d\'{e}riv%
\'{e}e de Lie $L_{X}$ le long d'un champ de vecteurs, section du fibr\'{e}
tangent $TM$) et de diff\'{e}rentielle ext\'{e}rieure $d_{\rho}$ (cf. \cite%
{Ana}).

Pour toute section $s$ du fibr\'{e} vectoriel $E$, il existe un unique
endomorphisme gradu\'{e} de degr\'{e} $0$ de l'alg\`{e}bre gradu\'{e}e
\underline{$\Lambda E^{\ast}$}, appel\'{e} d\'{e}riv\'{e}e de Lie
relativement \`{a} $s$ et not\'{e} $L_{s}^{\rho}$ v\'{e}rifiant les propri%
\'{e}t\'{e}s :

\begin{enumerate}
\item pour toute fonction $f\in\underline{\Lambda^{0}E^{\ast}}=\mathcal{F}$%
\begin{equation*}
L_{s}^{\rho}\left( f\right) =L_{\underline{\rho}\circ s}\left( f\right) =i_{%
\underline{\rho}\circ s}\left( df\right) \newline
\end{equation*}
o\`{u} $L_{X}$ d\'{e}signe la d\'{e}riv\'{e}e de Lie classique par rapport
au champ de vecteurs $X$

\item pour tout $q$--forme $\omega\in$\underline{$\Lambda^{q}E^{\ast}$} (o%
\`{u} $q>0$)%
\begin{align*}
\left( L_{s}^{\rho}\omega\right) \left( s_{1},\dots,s_{q}\right)
&=L_{s}^{\sigma}\left( \omega\left( s_{1},\dots,s_{q}\right) \right)\\
& -{%
\displaystyle\sum\limits_{i=1}^{q}} \omega\left( s_{1},\dots,s_{i-1},\left[
s,s_{i}\right] _{E},s_{i+1},\dots,s_{q}\right)
\end{align*}
\end{enumerate}

D'autre part, on d\'{e}finit aussi pour toute fonction $f\in\underline {%
\Lambda^{0}E^{\ast}}=\mathcal{F}$ l'\'{e}l\'{e}ment de \underline {$%
\Lambda^{1}E^{\ast}$}, not\'{e} $d_{\rho}f,$ par%
\begin{equation}
d_{\rho}f=\underline{t_{\rho}}\circ df  \label{d0}
\end{equation}

o\`{u} $t_{\rho}:T^{\ast}M\rightarrow E^{\ast}$ est la transpos\'{e}e de
l'ancre.

Il existe un unique endomorphisme gradu\'{e} de degr\'{e} $1$ de l'alg\`{e}%
bre gradu\'{e}e \underline{$\Lambda E^{\ast}$} appel\'{e} d\'{e}rivation ext%
\'{e}rieure, not\'{e} $d_{\rho}$, v\'{e}rifiant les propri\'{e}t\'{e}s~:

\begin{enumerate}
\item pour toute fonction $f\in\underline{\Lambda^{0}E^{\ast}}=\mathcal{F}$,
$d_{\rho}f$ est l'\'{e}l\'{e}ment de \underline{$\Lambda^{1}E^{\ast}$} d\'{e}%
fini \`{a} la relation (\ref{d0}).

\item pour tout \'{e}l\'{e}ment $\omega$ de \underline{$\Lambda^{q}E^{\ast}$}
($q>0$), $d_{\rho}\omega$ est l'unique \'{e}l\'{e}ment de \underline {$%
\Lambda^{q+1}E^{\ast}$} tel que pour tous $s_{0},\dots,s_{q}\in$\underline{$%
E $},
\begin{align*}
\left( d_{\rho}\omega\right) \left( s_{0},\dots,s_{q}\right) & ={%
\displaystyle\sum\limits_{i=0}^{q}}\left( -1\right) ^{i}L_{s_{i}}^{\rho
}\left( \omega\left( s_{0},\dots,\widehat{s_{i}},\dots,s_{q}\right) \right)
\\
& +{\displaystyle\sum\limits_{0\leq i<j\leq q}^{q}}\left( -1\right)
^{i+j}\left( \omega\left( \left[ s_{i},s_{j}\right] _{E},s_{0},\dots,%
\widehat{s_{i}},\dots,\widehat{s_{j}},\dots,s_{q}\right) \right)
\end{align*}
\end{enumerate}

On a alors la propri\'{e}t\'{e}
\begin{equation*}
d_{\rho}\circ d_{\rho}=0
\end{equation*}

\subsection{Morphismes d'alg\'{e}bro\"{\i}des}

\begin{definition}
Un morphisme de fibr\'{e}s vectoriels $\psi:E\rightarrow E^{\prime}$ au
dessus de $f:M\rightarrow M^{\prime}$ est un morphisme des alg\'{e}bro\"{\i}%
des de Lie banachiques $\left( E,\pi,M,\rho\right) $ et $\left(
E^{\prime},\pi^{\prime},M^{\prime},\rho^{\prime}\right) $ si l'application
$\psi^{\ast}:$\underline{$\Lambda^{q}E^{\prime\ast}$}$%
\rightarrow $\underline{$\Lambda^{q}E^{\ast}$} d\'{e}finie par
\begin{equation*}
\left( \psi^{\ast}\alpha^{\prime}\right) _{x}\left( s_{1},\dots,s_{q}\right)
=\alpha _{f\left( x\right) }^{\prime}\left( \psi\circ s_{1},\dots,\psi\circ
s_{q}\right)
\end{equation*}
commute avec les diff\'{e}rentielles~:%
\begin{equation*}
d_{\rho}\circ\psi^{\ast}=\psi\circ d_{\rho^{\prime}}
\end{equation*}
\end{definition}

\subsection{Courbes admissibles}

Dans le cadre de la M\'{e}canique, un \'{e}l\'{e}ment $a$ de $E$ peut \^{e}%
tre vu comme un vecteur vitesse g\'{e}n\'{e}ralis\'{e}, la vitesse naturelle
$v$ \'{e}tant obtenue par application de l'ancre $\rho$ \`{a} $a$, i.e. $%
v=\rho\left( a\right) $.

Une courbe $\gamma:\left[ 0,1\right] \rightarrow E$ est dite \textit{%
admissible} (cf. \cite{CLMM}) si $\overset{.}{m}\left( t\right) =\rho\left(
\gamma\left( t\right) \right) $ o\`{u} $t\mapsto m\left( t\right) =\pi\left(
\gamma\left( t\right) \right) $ est la courbe sur la base $M$.

Un morphisme d'alg\'{e}bro\"{\i}des de Lie applique alors les courbes
admissibles sur les courbes admissibles.

\subsection{Semi-gerbes}

Soit $\left(  E,\pi,M,\rho\right)  $ un alg\'{e}bro\"{\i}de de Lie banachique
et soit  $T\pi:TE\rightarrow TM$ l'application tangente de $\pi$. On notera
$\tau_{E}:TE\rightarrow E$ le fibr\'{e} tangent de $E$.

La notion de \textit{semi-gerbe} que l'on donne ici est une g\'{e}n\'{e}ralisation de
celle utilis\'{e}e lorsque $E=TM.$

\begin{definition}
Une section $S:E\rightarrow TE$ es appel\'{e}e une semi-gerbe si
\end{definition}

\begin{enumerate}
\item $\tau_{E}\circ S=\operatorname*{Id}_{E}$

\item $T\pi\circ S=\rho$
\end{enumerate}

Nous avons alors le lien suivant entre courbes admissibles et semi-gerbes (cf. \cite{Ana})

\begin{proposition}
Un champ de vecteurs sur $E$ est une semi-gerbe si et seulement si ses courbes
int\'{e}grables sont des courbes admissibles.
\end{proposition}

\section{\label{5_LimiteProjectiveForteAlgebroidesLieBanachiques}Limite
projective forte d'alg\'{e}bro\"{\i}des de Lie banachiques}

Un \textit{syst\`{e}me projectif fort d'alg\'{e}bro\"{\i}des de Lie
banachiques} est la donn\'{e}e d'une suite $\left(
E_{i},\pi_{i},M_{i},\rho_{i}\right) _{i\in\mathbb{N}}$ o\`{u}

\begin{description}
\item[--] $\left(  \left(  E_{i},f_{i}^{j}\right)  _{j\geq i}\right)
_{i\in\mathbb{N}}$ est un syst\`{e}me projectif fort de fibr\'{e}s vectoriels banachiques
($\pi_{i}:E_{i}\rightarrow M_{i}$) au-dessus du syst\`{e}me projectif fort  de vari\'{e}t\'{e}s
$\left(  \left(  M_{i},\delta_{i}^{j}\right)  _{j\geq i}\right)
_{i\in\mathbb{N}}$

\item[--] Pour tous $i,j\in\mathbb{N}$ tels que $j\geq i$, on a
\[
\rho_{i}\circ f_{i}^{j}=T\delta_{i}^{j}\circ\rho_{j}%
\]

\item[--] $f_{i}^{j}:E_{j}\rightarrow E_{i}$ est un morphisme d'alg\'{e}bro\"{\i}des de Lie $\left(  E_{j},\pi_{j},M_{j},\rho_{j}\right)  $ et $\left(
E_{i},\pi_{i},M_{i},\rho_{i}\right)  $
\end{description}

\begin{theorem}
\label{T_LimiteProjectiveAlgebroidesLieBanachiques}Soit $\left( E_{i},\pi
_{i},M_{i},\rho _{i}\right) _{i\in \mathbb{N}}$ un \textit{syst\`{e}me
projectif fort d'alg\'{e}bro\"{\i}des de Lie banachiques.}\newline
Alors $\left( \underleftarrow{\lim }E_{i},\underleftarrow{\lim }\pi _{i},%
\underleftarrow{\lim }M_{i},\underleftarrow{\lim }\rho _{i}\right) $ a une
structure d'alg\'{e}bro\"{\i}de de Lie fr\'{e}ch\'{e}tique.
\end{theorem}

\textbf{Preuve.---} Remarquons tout d'abord que la limite projective $\underleftarrow{\lim}M_{i}$
est munie d'une structure diff\'{e}rentiable comme d\'{e}finie en \ref{D_VarieteDifferentiableConvenient}. \\
$\left(  \lim E_{i},\underleftarrow{\lim}\pi_{i},\underleftarrow{\lim}M\right)  _{i}$ est un fibr\'{e} vectoriel fr\'{e}chetique de groupe structural $H^{0}\left(
\mathbb{E}\right)  $ (cf. Proposition \ref{P_LimProjFVB}). La limite projective des fibr\'{e}s (vetoriels) tangents $\left(  \underleftarrow
{\lim}TM_{i},\underleftarrow{\lim}p_{i},\underleftarrow{\lim}M_{i}\right)  $
est munie d'une structure de fibr\'{e} vectoriel fr\'{e}chetique~; on obtient alors le r\'{e}sultat de \cite{Gal3}, Theorem 2.1.

\bigskip

Etudions maintenant les propri\'{e}t\'{e}s des sections des fibr\'{e}s vectoriels $\underleftarrow{\lim}TM_{i}$,~ $\underleftarrow{\lim}E_{i}$ et la limite projective des ancres $\rho_{i}$.

Pour $\left(  g_{i}\right)  _{i\in\mathbb{N}}$ tels que $g_{j}=g_{i}%
\circ\delta_{i}^{j}=\left(  \delta_{i}^{j}\right)  ^{\ast}\left(
g_{i}\right)  $ on peut d\'{e}finir la limite projective $g=\underleftarrow{\lim
}g_{i}$ qui est encore de classe $C^\infty$.

Remarquons tout d'abord que si $X_{i}=T\delta_{i}^{j}\left(  X_{j}\right)  $, nous avons
$X_{i}\left(  g_{i}\right)  =\left(  T\delta_{i}^{j}\left(  X_{j}\right)
\right)  \left(  g_{i}\right)  =X_{j}\left(  g_{i}\circ\delta_{i}^{j}\right)
=X_{j}\left(  g_{j}\right)  $. On d\'{e}finit alors $X=\underleftarrow{\lim}X_{i}$
$\in\underleftarrow{\lim}\mathfrak{X}\left(  M_{i}\right)  $ et on obtient
$Xg=\underleftarrow{\lim}X_{i}g_{i}$ o\`{u} $X_{i}g_{i}\in\mathcal{F}_{i}$.

Si les suites $\left(  X_{i}^{1}\right)  _{i\in\mathbb{N}}$ et $\left(
X_{i}^{2}\right)  _{i\in\mathbb{N}}$ o\`{u} $X_{i}^{1},X_{i}^{2}\in
\mathfrak{X}\left(  M_{i}\right)  $ sont telles que $X_{i}^{1}=T\delta_{i}%
^{j}\left(  X_{j}^{1}\right)  $ (resp. $X_{i}^{2}=T\delta_{i}^{j}\left(
X_{j}^{2}\right)  $), elles donnenet naissance \`{a} des \'{e}l\'{e}ments $X^{1},X^{2}\in$
$\underleftarrow{\lim}\mathfrak{X}\left(  M_{i}\right)  $.

Puisque $X_{i}^{1}$ et $X_{j}^{1}$ sont $\delta_{i}^{j}-$reli\'{e}s (ainsi que
$X_{i}^{2}$ et $X_{j}^{2}$), leurs crochets $\delta_{i}^{j}-$ le sont aussi,
\newline i.e. $\left[  X_{i}^{1},X_{i}^{2}\right]  _{i}=T\delta_{i}^{j}\left(
\left[  X_{j}^{1},X_{j}^{2}\right]  _{j}\right)  $ et l'on obtient le crochet de
$X^{1}$ et $X^{2}$ comme limite projective de ces crochets.

\bigskip

Soit $s=\underleftarrow{\lim}s_{i}$ o\`{u} $s_{i}\in$\underline{$E_{i}$}$. $
Puisque les espaces $\underleftarrow{\lim}M_{i}$ et $\underleftarrow{\lim
}E_{i}$ sont des vari\'{e}t\'{e}s diff\'{e}rentiables, la section $s:\left(  x_{0},x_{1}%
,\dots\right)  \mapsto\left(  s_{0}\left(  x_{0}\right)  ,s_{1}\left(
x_{1}\right)  ,\dots\right)  $ est de classe $C^\infty$ (cf. definition
\ref{D_AppClasseCInfini}).

\bigskip

Prouvons maintenant que l'on peut d\'{e}duire la condition de compatibilit\'{e}%
\[
f_{i}^{j}\circ\left[  s_{j}^{1},s_{j}^{2}\right]  _{E_{j}}=\left[  s_{i}%
^{1},s_{i}^{2}\right]  _{E_{i}}\circ\delta_{i}^{j}%
\]

de la structure de morphisme de  $f_{i}^{j}$ (commutativit\'{e} avec les diff\'{e}rentielles appliqu\'{e}e aux $1-$formes)

Nous avons $\left(  \left(  f_{i}^{j}\right)  ^{\ast}\left(  d_{E_{i}}\alpha
_{i}\right)  \right)  \left(  s_{j}^{1},s_{j}^{2}\right)  =\left(  d_{E_{i}%
}\alpha_{i}\right)  \left(  f_{i}^{j}\circ s_{j}^{1},f_{i}^{j}\circ s_{j}%
^{2}\right)  $ o\`{u}%
\begin{align*}
&  \left(  d_{E_{i}}\alpha_{i}\right)  \left(  f_{i}^{j}\circ s_{j}^{1}%
,f_{i}^{j}\circ s_{j}^{2}\right) \\
&  =L_{\rho_{i}\circ\left(  f_{i}^{j}\circ s_{j}^{1}\right)  }\left(
\alpha_{i}\left(  f_{i}^{j}\circ s_{j}^{2}\right)  \right)  -L_{\rho_{i}%
\circ\left(  f_{i}^{j}\circ s_{j}^{2}\right)  }\left(  \alpha_{i}\left(
f_{i}^{j}\circ s_{j}^{1}\right)  \right)  -\alpha_{i}\left[  f_{i}^{j}\circ
s_{j}^{1},f_{i}^{j}\circ s_{j}^{2}\right]  _{E_{i}}\\
&  =L_{\rho_{i}\circ s_{i}^{1}}\left(  \alpha_{i}\left(  s_{i}^{2}\right)
\right)  -L_{\rho_{i}\circ s_{i}^{2}}\left(  \alpha_{i}\left(  s_{i}%
^{1}\right)  \right)  -\alpha_{i}\left[  s_{i}^{1},s_{i}^{2}\right]  _{E_{i}%
}\\
&  =X_{i}^{1}\left(  \alpha_{i}\left(  s_{i}^{2}\right)  \right)  -X_{i}%
^{2}\left(  \alpha_{i}\left(  s_{i}^{1}\right)  \right)  -\alpha_{i} \left[
s_{i}^{1},s_{i}^{2}\right]  _{E_{i}}%
\end{align*}

D'un autre c\^{o}t\'{e}, nous avons%

\begin{align*}
&  \left(  d_{E_{j}}\left(  \left(  f_{i}^{j}\right)  ^{\ast}\alpha
_{i}\right)  \right)  \left(  s_{j}^{1},s_{j}^{2}\right) \\
&  =L_{\rho_{j}\circ s_{j}^{1}}\left(  \left(  \left(  f_{i}^{j}\right)
^{\ast}\alpha_{i}\right)  \left(  s_{j}^{2}\right)  \right)  -L_{\rho_{j}\circ
s_{j}^{2}}\left(  \left(  \left(  f_{i}^{j}\right)  ^{\ast}\alpha_{i}\right)
\left(  s_{j}^{1}\right)  \right)  -\left(  \left(  f_{i}^{j}\right)  ^{\ast
}\alpha_{i}\right)  \left[  s_{j}^{1},s_{j}^{2}\right]  _{E_{j}}\\
&  =X_{j}^{1}\left(  \alpha_{i}\left(  f_{i}^{j}\circ s_{j}^{2}\right)
\right)  -X_{j}^{2}\left(  \alpha_{i}\left(  f_{i}^{j}\circ s_{i}^{1}\right)
\right)  -\alpha_{i}\left[  f_{i}^{j}\circ s_{j}^{1},f_{i}^{j}\circ s_{i}%
^{2}\right]  _{E_{i}}\\
&  =X_{i}^{1}\left(  \alpha_{i}\left(  s_{i}^{2}\right)  \right)  -X_{i}%
^{2}\left(  \alpha_{i}\left(  s_{i}^{1}\right)  \right)  -\alpha_{i} \left[
f_{i}^{j}\circ s_{j}^{1},f_{i}^{j}\circ s_{i}^{2}\right]  _{E_{i}}%
\end{align*}

Puisque $f_{i}^{j}$ est un morphisme, on obtient $\alpha_{i}\left[  s_{i}^{1}%
,s_{i}^{2}\right]  _{E_{i}}=\alpha_{i}\left[  f_{i}^{j}\circ s_{j}^{1}%
,f_{i}^{j}\circ s_{i}^{2}\right]  _{E_{i}}$ et la condition de compatibilit\'{e}.

Ainsi le crochet  $\left[  s^{1},s^{2}\right]  _{\underleftarrow{\lim}E_{i}}$ des limites projectives de sections $s^{1}=\underleftarrow{\lim}s_{i}^{1}$ et
$s^{2}=\underleftarrow{\lim}s_{i}^{2}$ peut-\^{e}tre d\'{e}fini comme limite projective des sections $\left[  s_{i}^{1},s_{i}^{2}\right]  _{E_{i}}$ de $E_{i}$.

L'ensemble $\underleftarrow{\lim}E_{i}$ muni de ce crochet a alors une structure d'alg\`{e}bre de Lie.

Eu \'{e}gard aux conditions $\rho_{i}\circ f_{i}^{j}=T\delta_{i}^{j}\circ
\rho_{j}$ la limite projective $\rho=\underleftarrow{\lim}\rho_{i}$ est un morphisme de fibr\'{e}s.

L'application $\rho=\underleftarrow{\lim}\rho_{i}$ est alors un morphisme d'alg\`{e}bres de Lie entre $\left(  \underleftarrow{\lim}\underline{E_{i}},\left[
.,.\right]  _{\underleftarrow{\lim}E_{i}}\right)  $ et $\left(
\underleftarrow{\lim}\underline{TM_{i}},\left[  .,.\right]  _{i}\right)  .$

Pour tous $i\in\mathbb{N}$, chaque section $s_{i}^{1}$ et $s_{i}^{2}$ of $E_{i}$
et chaque fonction de classe $C^\infty$ $g_{i}:M_{i}\rightarrow\mathbb{R},$ on a%
\[
\left[  s_{i}^{1},g_{i}s_{i}^{2}\right]  _{E_{i}}=g_{i}\left[  s_{i}^{1}%
,s_{i}^{2}\right]  _{E_{i}}+\left(  \rho_{i}\left(  s_{i}^{1}\right)  \right)
\left(  g_{i}\right)  \ s_{i}^{2}%
\]

Afin d'obtenir la relation:%
\[
\left[  s_{1},gs_{2}\right]  _{E}=g\left[  s_{1},s_{2}\right]  _{E}+\left(
\rho\left(  s_{1}\right)  \right)  \left(  g\right)  \ s_{2}%
\]

nous devons prouver que~:

1) $f_{i}^{j}\circ\left(  g_{j}\left[  s_{j}^{1},s_{j}^{2}\right]  \right)
=g_{i}\left[  s_{i}^{1},s_{i}^{2}\right]  \circ\delta_{i}^{j}$

2) $f_{i}^{j}\circ\left[  \left(  \rho_{j}\left(  s_{j}^{1}\right)  \right)
\left(  g_{j}\right)  \ s_{j}^{2}\right]  =\left[  \left(  \rho_{i}\left(
s_{i}^{1}\right)  \right)  \left(  g_{i}\right)  \ s_{i}^{2}\right]
\circ\delta_{i}^{j}$

En ce qui concerne le premier point, pour chaque fil $\left(  x_{i}\right)  _{i\in N}$ , i.e.
$x_{j}=\delta_{i}^{j}\left(  x_{j}\right)  ,$ nous avons%
\[
f_{i}^{j}\circ\left(  g_{j}\left[  s_{j}^{1},s_{j}^{2}\right]  _{E_{j}%
}\right)  \left(  x_{j}\right)  =f_{i}^{j}\left(  \left(  g_{i}\circ\delta
_{i}^{j}\times\left[  s_{j}^{1},s_{j}^{2}\right]  _{E_{j}}\right)  \left(
x_{j}\right)  \right)
\]

Puisque $f_{i}^{j}$ est une application lin\'{e}aire de $\pi_{j}^{-1}\left(
x_{j}\right)  $ vers $\pi_{i}^{-1}\left(  x_{i}\right)  ,$ cette expression
est \'{e}gale \`{a} $g_{i}\left(  x_{i}\right)  \times f_{i}^{j}\left(  \left[  s_{j}%
^{1},s_{j}^{2}\right]  _{E_{j}}\left(  x_{j}\right)  \right)  $. Gr\^{a}ce \`{a} la condition de compatibilit\'{e} $f_{i}^{j}\circ\left[  s_{j}^{1},s_{j}^{2}\right]
_{E_{j}}=\left[  s_{i}^{1},s_{i}^{2}\right]  _{E_{i}}\circ\delta_{i}^{j}$ nous avons \'{e}tabli le premier point.

Pour ce qui est du second point, nous utilisons tout d'abord la commutativit\'{e} avec les diff\'{e}rentielles
$d_{E_{i}}$ and $d_{E_{j}}$.\\
\[
\left[  \left(  f_{i}^{j}\right)  ^{\ast}\left(  d_{E_{i}}g_{i}\right)
\right]  \left(  s_{j}\right)  \left(  x_{j}\right)  =\left[  d_{E_{j}}\left(
\left(  \delta_{i}^{j}\right)  ^{\ast}\left(  g_{i}\right)  \right)  \right]
\left(  s_{j}\right)  \left(  x_{j}\right)
\]

et ainsi

\[
\left(  d_{E_{i}}g_{i}\right)  \left(  f_{i}^{j}\circ s_{j}\right)  \left(
x_{j}\right)  =\left[  d_{E_{j}}\left(  g_{j}\right)  \right]  \left(
s_{j}\right)  \left(  x_{j}\right)
\]

En utilisant la d\'{e}finition $d_{E_{i}},$ i.e. $d_{E_{i}}g_{i}=$ $t_{\rho_{i}}\circ
dg_{i}$ on obtient%
\[
dg_{i}\left[  \rho_{i}\left(  f_{i}^{j}\left(  s_{j}\left(  x_{j}\right)
\right)  \right)  \right]  =dg_{j}\left(  \rho_{j}\left(  s_{j}\left(
x_{j}\right)  \right)  \right)
\]

et donc%
\[
\left[  \rho_{i}\left(  f_{i}^{j}\circ s_{j}\right)  \right]  \left(
g_{i}\right)  \left(  x_{i}\right)  =\left[  \rho_{j}\circ s_{j}\right]
\left(  g_{j}\right)  \left(  x_{j}\right)
\]

Gr\^{a}ce \`{a} la condition de compatibilit\'{e}, nous obtenons
\[
\left[  f_{i}^{j}\circ\left(  \rho_{j}\circ s_{j}\right)  \right]  \left(
g_{i}\right)  =\left[  \rho_{j}\left(  s_{j}\right)  \right]  \left(  \left(
\delta_{i}^{j}\right)  ^{\ast}\circ g_{i}\right)
\]

Le second point en d\'{e}coule facilement. \hfill Q.E.D.

\section{\label{6_Exemples}Exemples}

\subsection{Alg\'{e}bro\"{\i}des de Nijenhuis Lie }

Soit~$\left(  \left(  M_{i},\delta_{i}^{j}\right)  _{j\geq i}\right)
_{i\in\mathbb{N}}$ un syst\`{e}me projectif fort de vari\'{e}t\'{e}s banachiques. Pour tout
$i\in\mathbb{N}$, consid\'{e}rons un tenseur de Nijenhuis  $N_{i}:TM_{i}\rightarrow TM_{i}
$ (cf. exemple \ref{EX_AlgebroideLie_TenseurNijenhuis}). Dans ce cas, consid\'{e}rons $f_{i}^{j}=T\delta_{i}^{j}$ morphisme de $TM_{j}$ dans $TM_{i}$. Si nous avons la condition de compatibilit\'{e}%
\[
N_{i}\circ T\delta_{i}^{j}=T\delta_{i}^{j}\circ N_{j}%
\]

$\left(  \underleftarrow{\lim}TM_{i},\underleftarrow{\lim}\pi_{i}%
,\underleftarrow{\lim}M_{i},\underleftarrow{\lim}N_{i}\right)  $ est un alg\'{e}bro\"{\i}de de Lie
fr\'{e}chetique Lie puisque nous avons en particulier%
\[
\left(  \delta_{i}^{j}\right)  ^{\ast}\circ d_{TM_{i}}=d_{TM_{j}}\circ\left(
\delta_{i}^{j}\right)  ^{\ast}.
\]

\bigskip

A titre d'exemple, on peut considérer le case de l'oscillateur harmonique de dimension infinie qui est un syst\`{e}me hamiltonien $L-$integrable au sens de \cite{Liu}. On consid\`{e}re la limite projective $\left(  \left(  \mathbb{R}^{2i},\delta_{i}^{j}\right)  _{j\geq i}\right)  _{i\in\mathbb{N}}$ o\`{u} $\delta_{i}^{j}$ est la projection canonique de  $\mathbb{R}^{2j}$ sur $\mathbb{R}^{2i}$. Le tenseur de Nijenhuis  $N_{i}$ qui correspond \`{a} l'op\'{e}rateur de r\'{e}cursion peut \^{e}tre \'{e}crit comme suit~:%
\[
N_{i}=
{\displaystyle\sum\limits_{k=1}^{i}}
\left(  x_{k}^{2}+y_{k}^{2}\right)  \left(  dx_{k}\otimes\dfrac{\partial
}{\partial x_{k}}+dy_{k}\otimes\dfrac{\partial}{\partial y_{k}}\right)
\]

o\`{u} $\left(  \left(  x_{1},y_{1}\right)  ,\dots,\left(  x_{i},y_{i}\right)
\right)  $ sont les coordonn\'{e}es sur  $\mathbb{R}^{2i}$. Il est alors facile d'\'{e}tablir la condition de compatibilit\'{e}.

\subsection{Distributions}

\subsubsection{Limite projective de distributions involutives}

 Une distribution sur une vari\'{e}t\'{e} banachique $B$ est une application $C^\infty$  $D:B \rightarrow TB$ telle que pour tout $x\in B, D_x$ est un sous espace vectoriel $T_{x}B$. Cette distribution est involutive si pour tous champs de vecteurs $X$ et $Y$ tangents \`{a} $D$, le crochet $\left[
X,Y\right]$ est encore tangent \`{a} $D$.

On peut remarquer que la distribution image de l'ancre d'un alg\'{e}bro\"{\i}de de Lie est une distribution involutive faible (cf. \cite{Pel}) si la base est r\'{e}guli\`{e}ment lisse (cf. \cite{CabPel}).\\

Soit~$\left(  \left(  M_{i},\delta_{i}^{j}\right)  _{j\geq i}\right)
_{i\in\mathbb{N}}$ un syst\`{e}me projectif fort de vari\'{e}t\'{e}s banachiques.
Consid\'{e}rons pour tout $i\in\mathbb{N}$, une distribution involutive diff\'{e}rentiable $E_{i}$ au dessus de la vari\'{e}t\'{e} $M_{i}$.

$\left(  E_{i},\pi_{i},M_{i},J_{i}\right)  _{i\in\mathbb{N}}$, o\`{u}
$J_{i}:E_{i} \rightarrow TM_{i}$ est l'injection naturelle et o\`{u} $f_{i}^{j}$ est la restriction de $T\delta_{i}^{j}$ \`{a} $E_{i}$, est un syst\`{e}me projectif fort d'alg\'{e}bro\"{\i}des de Lie.

La limite projective $\underleftarrow{\lim}E_{i}$ peut \^{e}tre vue comme une distribution involutive sur le fibr\'{e} fr\'{e}ch\'{e}tique  $\underleftarrow{\lim}TM_{i}$.\\

\subsubsection{Limite projective de distributions de rangs finis}
Consid\'{e}rons  le cas d'une distribution de dimension $1$ sur l'ensemble des jets infinis de sections d'un fibr\'{e} vectoriel $p:F \rightarrow N$. Soit $X$ un champ de vecteurs sur $F$ projetable sur $N$ de projection $\hat{X}$; le flot $\varphi_{t}^{X}$ de $X$ est au dessus du flot $\varphi_{t}^{\hat{X}}$ de $\hat{X}$ et par prolongement \`{a} $J^\infty (F)$ on obtient un groupe local \`{a} un param\`{e}tre $\phi_t=pr^\infty(\varphi_{t}^{X})$ de transformations sur $J^\infty (F)$ (cf. \cite{And}, \cite{Olv}). Le prolongement $pr^\infty(X)$ du champ de vecteurs $X$ est le champ de vecteurs sur $J^\infty(p)$ associ\'{e} \`{a} ce flot. Qui plus est ce flot pr\'{e}serve la distribution de Cartan (id\'{e}al de contact) $\cal C$.

On peut remarquer que $\cal C$ est une  distribution involutive sur la limite projective
$\underleftarrow{\lim}TJ^{i}(p)$ qui appara\^{\i}t comme limite de distributions non involutives $J^{i}(p)$
(cf. \cite{Sau}).\\
Rappelons qu'une vari\'{e}té\'{e} diff\'{e}rentiable fr\'{e}ch\'{e}tique munie d'une distribution involutive correspond \`{a} la notion de  diffi\'{e}t\'{e} (eng. diffiety: differential variety) comme introduit par Vinogradov dans (\cite{Vin}).
On peut trouver des applications d'une telle situation en m\'{e}canique non holonome et en th\'{e}orie du contr\^{o}le non lin\'{e}aire (cf. par exemple \cite{FLMR}).

\bigskip

Si l'on consid\`{e}re un syst\`{e}me d'EDP $\mathcal{E}$, i.e. une sous vari\'{e}t\'{e} du fibr\'{e} des jets $J^k(\pi)$, on obtient, par prolongement d'ordre infini, une sous vari\'{e}t\'{e}  $i:\mathcal{E} \rightarrow J^{\infty}\left(\pi\right)$ de $J^{\infty}\left(\pi\right)$. On a alors une distribution involutive sur $\mathcal{E}$ par restriction de la distribution de Cartan \`{a} $\mathcal{E}$ via l'image r\'{e}ciproque par $i$ (cf. \cite{KisVan}, \cite{Dri}).\\

\bigskip

Consid\'{e}rons maintenant le cas d'un syst\`{e}me d'EDP d'ordre $1$ \`{a} $2$
variables. Soit le fibr\'{e} trivial $\pi:E\rightarrow M$ o\`{u}
$M=\mathbb{R}^{2}$ et $E=M\times\mathbb{R}^{2}=\mathbb{R}\times\mathbb{R}^{2}%
$. On munit le fibr\'{e} $J^{1}\left(  \pi\right)  $ des jets d'ordre $1$ de
sections lisses de $E$ au dessus de $M$ du syst\`{e}me de coordonn\'{e}es
$\left(  x,y,u,u_{x},u_{y}\right)  $, les variables ind\'{e}pendantes sont $x$
et $y$ et la variable d\'{e}pendante est $u$.\newline On consid\`{e}re alors
le syst\`{e}me d'EDP du premier ordre%
\[
\mathcal{R}_{1}\quad\left\{
\begin{array}
[c]{c}%
u_{x}=\varphi\left(  x,y,u\right)  \\
u_{y}=\psi\left(  x,y,u\right)
\end{array}
\right.
\]

Une condition n\'{e}cessaire pour que ce syst\`{e}me soit int\'{e}grable est
que les conditions d'int\'{e}grabilit\'{e} suivantes soit v\'{e}rifi\'{e}es
sur la sous-vari\'{e}t\'{e} $\mathcal{R}_{1}$ $:$%
\[
D_{x}\psi-D_{y}\varphi\ _{|\mathcal{R}_{1}}=0
\]

Dans le cadre d'un syst\`{e}me de ce type, on \'{e}tablit que ces conditions
sont aussi suffisantes (cf. \cite{Sei}, 1.4) .

Ces conditions d'int\'{e}grabilit\'{e} s'\'{e}crivent alors%
\begin{equation}
\dfrac{\partial\psi}{\partial x}+\varphi\dfrac{\partial\psi}{\partial
u}-\dfrac{\partial\varphi}{\partial y}-\psi\dfrac{\partial\varphi}{\partial
u}=0\tag{CI}%
\end{equation}

Ces conditions d'int\'{e}grabilit\'{e} peuvent \^{e}tre interpr\'{e}t\'{e}es
de mani\`{e}re plus g\'{e}om\'{e}trique en introduisant les deux champs de
vecteurs sur $E$%
\[
X=\dfrac{\partial}{\partial x}+\varphi\dfrac{\partial}{\partial u}\text{ et
}Y=\dfrac{\partial}{\partial y}+\psi\dfrac{\partial}{\partial u}%
\]
\newline Ces deux champs engendrent une distribution r\'{e}guli\`{e}re
$\mathcal{D}^{0}$ de dimension $2$.\newline La condition
d'int\'{e}grabilit\'{e} (CI) est alors \'{e}quivalente \`{a} $\left[
X,Y\right]  =0$, i.e. \`{a} l'involutivit\'{e} de la distribution.\\

Le processus de prolongement (cf. \cite{Olv}, 5.1.) permet de prolonger chacun
des champs de vecteurs $X$ (resp. $Y$)  sur $J^{1}\left(  \pi\right)  ,$
$J^{2}\left(  \pi\right)  $ ainsi que $J^{\infty}\left(  \pi\right)  $ en
utilisant les op\'{e}rateurs de d\'{e}rivation totale (cf. \cite{Olv}, 2.3.)
en les champs $\operatorname*{pr}\nolimits^{\left(  1\right)  }X$,
$\operatorname*{pr}\nolimits^{\left(  2\right)  }X,\dots,\operatorname*{pr}%
\nolimits^{\infty}X$ (resp. $\operatorname*{pr}\nolimits^{\left(  1\right)
}Y$, $\operatorname*{pr}\nolimits^{\left(  2\right)  }Y,\dots
,\operatorname*{pr}\nolimits^{\infty}Y$)

Eu \'{e}gard \`{a} la formule g\'{e}n\'{e}rale \cite{Olv}, (2.47)%
\[
\operatorname*{pr}\nolimits^{\left(  n\right)  }\left[  X,Y\right]  =\left[
\operatorname*{pr}\nolimits^{\left(  n\right)  }X,\operatorname*{pr}%
\nolimits^{\left(  n\right)  }Y\right]
\]
\newline\newline on en d\'{e}duit que la distribution $\mathcal{D}^{1}$
engendr\'{e}e par les premiers prolongements $\operatorname*{pr}%
\nolimits^{\left(  1\right)  }X$ et $\operatorname*{pr}\nolimits^{\left(
1\right)  }Y$ est encore involutive, de m\^{e}me que les diverses
distributions $\mathcal{D}^{n}$ engendr\'{e}e par les prolongements d'ordre
$n$ que sont $\operatorname*{pr}\nolimits^{\left(  n\right)  }X$ et
$\operatorname*{pr}\nolimits^{\left(  n\right)  }Y$. La limite projective est
alors une diffi\'{e}t\'{e}.

\subsection{Limite inverse de distributions banachiques}
On consid\`{e}re ici le cas o\`{u} les applications $\delta_{i}^{i+1}:M_{i+1} \rightarrow M_i$ sont les injections canoniques entre vari\'{e}t\'{e}s de Banach, les distributions $E_i$ de corang $1$ sont d\'{e}finies comme $\ker \alpha_i$ o\`{u} $\alpha_i$ est une $1$-forme v\'{e}rifiant les diff\'{e}rentes conditions de compatibilit\'{e}. \\
On peut rencontre cette situation pour $M_i=C^{i}(\mathbb{S}^1)$ o\`{u} $\alpha_i(u_i)=\int_{\mathbb{S}^1}u_i(x)dx$. La distribution associ\'{e}e est affine (\cite{KapMak}) et est li\'{e}e au premier tenseur de Poisson de l'\'{e}quation de KdV.

\section{\label{7_LimiteProjectiveSemisprays}Limite projective forte de semi-gerbes}

Soit $\left(  E_{i},\pi_{i},M_{i},\rho_{i}\right)  _{i\in\mathbb{N}}$ un syst\`{e}me projectif fort d'alg\'{e}bro\"{\i}des de Lie.

%%%%%%%%%%%%%%%%%  SPL of admissible curves
Consid\'{e}rons une suite $\left(  \gamma_{i}\right)  _{i\in\mathbb{N}}$ o\`{u}
$\gamma_{i}:\left[  0,1\right]  \rightarrow E_{i}$ est une courbe admissible telle que pour tous $i,j\in\mathbb{N}$ tels que $j\geq i$
\[
f_{i}^{j}\circ\gamma_{j}=\gamma_{i}%
\]

Par cons\'{e}quent $\gamma=\underleftarrow{\lim}\gamma_{i}$ existe.

Pour tous $i,j\in\mathbb{N}$ tels que $j\geq i$ et pour tous $t\in\left[
0,1\right]  ,$ en faisant appel aux \'{e}galit\'{e}s
\[
\left(  \pi_{i}\circ\gamma_{i}\right)  ^{\prime}\left(  t\right)  =\left(
\delta_{i}^{j}\circ\left(  \pi_{j}\circ\gamma_{j}\right)  \right)  ^{\prime
}\left(  t\right)  =T\delta_{i}^{j}\left(  \left(  \pi_{j}\circ\gamma
_{j}\right)  ^{\prime}\left(  t\right)  \right)
\]

on obtient :
\[
\left(  \pi_{i}\circ\gamma_{i}\right)  ^{\prime}\left(  t\right)  -\rho
_{i}\left(  \gamma_{i}\left(  t\right)  \right)  =T\delta_{i}^{j}\left(
\left(  \pi_{j}\circ\gamma_{j}\right)  ^{\prime}\left(  t\right)  -\rho
_{j}\left(  \gamma_{j}\left(  t\right)  \right)  \right)
\]

Et ainsi pour tout $t\in\left[  0,1\right]  $, nous avons :
\[
\left(  \pi\circ\gamma\right)\prime  \left(  t\right)  =\rho\left(  \gamma\left(
t\right)  \right)
\]

Une telle courbe sera appel\'{e}e \textit{courbe admissible} dans $E=$ $\underleftarrow
{\lim}E_{i}$.

\bigskip
%%%%%%%%%%%%%%%%%  LPF de semi-gerbes
Consid\'{e}rons maintenant $\left(  S_{i}\right)  _{i\in\mathbb{N}}$ o\`{u}
$S_{i}:E_{i}\rightarrow TE_{i}$ est une semi-gerbe tel que
\[
Tf_{i}^{j}\circ S_{j}=S_{i}\circ f_{i}^{j}%
\]

On peut alors d\'{e}finir $S:\left(  u_{0},u_{1},\dots\right)  \mapsto\left(
S_{0}\left(  u_{0}\right)  ,S_{1}\left(  u_{1}\right)  ,\dots\right)  $ qui est une section $C^\infty$ de   $\underleftarrow{\lim}TE_{i}$.

Il est facile de voir que l'on a  $\tau_{E}\circ S=Id_{E}$. \\

Pour tous $i,j\in\mathbb{N}$ tels que $j\geq i$ et tout $u_{i}=f_{i}^{j}\left(  u_{j}\right)  $ on a :

\[%
\begin{array}
[c]{ll}%
\left(  T\pi_{i}\circ s_{i}-\rho_{i}\right)  \left(  u_{i}\right)   & =\left(
T\pi_{i}\circ s_{i}\circ f_{i}^{j}-\rho_{i}\circ f_{i}^{j}\right)  \left(
u_{j}\right)  \\
& =\left(  T\pi_{i}\circ Tf_{i}^{j}\circ s_{j}-\rho_{i}\circ f_{i}^{j}\right)
\left(  u_{j}\right)  \\
& =\left(  T\left(  \pi_{i}\circ f_{i}^{j}\right)  \circ s_{j}-\rho_{i}\circ
f_{i}^{j}\right)  \left(  u_{j}\right)  \\
& =\left(  T\left(  \delta_{i}^{j}\circ\pi_{j}\right)  \circ s_{j}-\rho
_{i}\circ f_{i}^{j}\right)  \left(  u_{j}\right)
\end{array}
\]

Finalement, on peut \'{e}crire
\[
\left(  T\pi_{i}\circ s_{i}-\rho_{i}\right)  \left(  u_{i}\right)
=T\delta_{i}^{j}\left(  T\pi_{j}\circ s_{j}-\rho_{j}\right)  \left(
u_{j}\right)
\]

Et ainsi nous obtenons
\[
T\pi\circ s=\rho
\]

$S$ sera appel\'{e}e \textit{semi-gerbe}.
\bigskip

On peut d\'{e}finir de mani\`{e}re \'{e}vidente la notion de \textit{gerbe} sur la limite projective $\underleftarrow{\lim}E_{i}$ comme limite projective de gerbes..

\bigskip
%%%%%%%%%%%%%%%%%%  lien semi-gerbe - admissible curve

On termine ce travail par une proposition qui fait appara\^{\i}tre le lien entre semi-gerbes et courbes admissibles. On g\'{e}n\'{e}ralise alors un résulatat de  \cite{Ana}  (pour le cas des gerbes, dans le cas particulier $E=TM$, voir \cite{RezMal}).\\

\begin{proposition}
Un champ de vecteurs $S=\underleftarrow{\lim}S_{i}$ on $E=\underleftarrow{\lim}E_{i}$ est une semi-gerbe si et seulement si toutes ses courbes int\'{e}grales sont des courbes admissibles.
\end{proposition}

\textbf{Preuve.---} La preuve n'est rien d'autre qu'une adaptation de la d\'{e}monstration du Th\'{e}or\`{e}me 2.3 que l'on peut trouver dans l'article \cite{Ana}.
Consid\'{e}rons pour ce faire une semi-gerbe $S=\underleftarrow{\lim}S_{i}$ et supposons que $c:\left[  0,1\right]  \rightarrow E$ est une courbe int\'{e}grale de $S$.
Alors pour tout $i\in\mathbb{N},$ $c_{i}:\left[  0,1\right]  \rightarrow E_{i}$ est une courbe int\'{e}grale de $S_{i}$ (i.e. $\forall t\in\left[  0,1\right]  ,c_{i}^{\prime}\left(
t\right)  =S_{i}\left(  c_{i}\left(  t\right)  \right)  $) o\`{u} $f_{i}%
^{j}\circ c_{j}=c_{i}.$ Il s'ensuit que pour tout $i\in\mathbb{N}$ et pour tout
$t\in\left[  0,1\right]  $ on a :
\[
T\pi_{i}\circ c_{i}^{\prime}\left(  t\right)  =\left(  T\pi_{i}\circ
S_{i}\right)  \left(  c_{i}\left(  t\right)  \right)
\]
Puisque $\pi_{i}\circ c^{\prime}\left(  t\right)  =\rho_{i}\left(  c\left(
t\right)  \right)  $, $c_{i}$ est une courbe admissible, il en est de m\^{e}me pour
$c=\underleftarrow{\lim}c_{i}$. \\
La r\'{e}ciproque est laiss\'{e}e au lecteur. \hfill  Q.E.D.

\bigskip

Pour un syst\`{e}me projectif de gerbes il est facile de prouver, en utilisant la relation $h_{i}^{\lambda}\circ f_{i}^{j}=f_{i}^{j}\circ h_{j}^{\lambda}$ et les propri\'{e}t\'{e}s des applications tangentes, que pour tous $i,j\in\mathbb{N}$ tels que $j\geq i$
et pour tous  $u_{i}=f_{i}^{j}\left(  u_{j}\right) $ on a :

\[
Tf_{i}^{j}\left(  s_{j}\circ h_{j}^{\lambda}-\lambda Th_{j}^{\lambda}\circ
s_{j}\right)  \left(  u_{j}\right)  =\left(  s_{i}\circ h_{i}^{\lambda
}-\lambda Th_{i}^{\lambda}\circ s_{i}\right)  \left(  u_{i}\right)
\]
 On peut alors \'{e}crire $S\circ h_{\lambda}=\lambda Th_{\lambda}\circ S$ et $S$ est alors une gerbe sur $E$.

\end{document}